\documentclass[letterpaper, 10 pt, conference]{ieeeconf}  

\IEEEoverridecommandlockouts                              

\usepackage{graphics} 
\usepackage{epsfig} 
\usepackage{mathptmx} 
\usepackage{times} 
\usepackage{amsmath} 
\usepackage{amssymb}  
\usepackage{amsbsy}
\usepackage{dsfont}
\usepackage{algorithm}
\usepackage{algpseudocode}
\usepackage{url}
\graphicspath{{./images/}}

\DeclareMathOperator*{\argmin}{arg\,min}
\newcommand{\R}{\mathbb{R}}
\newcommand{\B}[1]{\mathbf{#1}}
\newcommand{\Ba}{\boldsymbol{a}}
\newcommand{\1}{\mathds{1}}
\newcommand{\innerprod}[2]{\left\langle #1,#2 \right\rangle}
\title{\LARGE \bf
An Efficient Semi-Real-Time Algorithm for Path Planning in the Hamilton-Jacobi Formulation
}

\author{Christian Parkinson$^{1}$ and Kyle Polage$^{2}$
\thanks{*This work was supported by the NSF through the Research Training Group on Applied Mathematics and Statistics for Data Driven Discovery at the University of Arizona (DMS-1937229).}
\thanks{$^{1}$Christian Parkinson is a postdoctoral associate with the Department of Mathematics, University of Arizona, Tucson, AZ 85721
        {\tt\small chparkin@math.arizona.edu}}%
\thanks{$^{2}$Kyle Polage is a student in the Department of Mathematics, Washington State University, Pullman, WA 99163
        {\tt\small kyle.polage@wsu.edu}}%
}

\begin{document}

\maketitle
\thispagestyle{empty}
\pagestyle{empty}

\begin{abstract}
We present a semi-real-time algorithm for minimal-time optimal path planning based on optimal control theory, dynamic programming, and Hamilton-Jacobi (HJ) equations. Partial differential equation (PDE) based optimal path planning methods are well-established in the literature, and provide an interpretable alternative to black-box machine learning algorithms. However, due to the computational burden of grid-based PDE solvers, many previous methods do not scale well to high dimensional problems and are not applicable in real-time scenarios even for low dimensional problems. We present a semi-real-time algorithm for optimal path planning in the HJ formulation, using grid-free numerical methods based on Hopf-Lax formulas. In doing so, we retain the intepretablity of PDE based path planning, but because the numerical method is grid-free, it is efficient and does not suffer from the curse of dimensionality, and thus can be applied in semi-real-time and account for realistic concerns like obstacle discovery. This represents a significant step in averting the tradeoff between interpretability and efficiency. We present the algorithm with application to synthetic examples of isotropic motion planning in two-dimensions, though with slight adjustments, it could be applied to many other problems. 
\end{abstract}

\section{INTRODUCTION}

With the proliferation of unmanned vehicles, automated navigation, and many other applications in robotics, the problem of optimal trajectory generation has become increasingly important. Especially in high-leverage applications like self-driving cars, it is vital to develop motion-planning methods which are efficient, accurate, and interpretable. 

Many state-of-the-art path planning methods rely heavily on deep learning \cite{DeepLearning1,DeepLearning2,DeepLearning5}. Deep learning is a powerful tool for such problems due to efficiency and flexibility with respect to complex constraints and modeling concerns. However, as is now well-documented, deep neural nets can suffer from lack of robustness and interpretability \cite{Hamon}, which presents an issue in decision critical problems. While robust, interpretable machine learning is an active area of research \cite{TowardRobustInterpretable,Ross,PIML}, routing models which do not rely on learning architectures provide an interpretable alternative.

One such collection of methods is partial differential equation (PDE), optimal control, and dynamic programming based planning algorithms. At a basic level, many such methods can be seen as continuous extensions of Dijkstra's famous algorithm \cite{Dijk} for traversing a weighted graph. These include level-set methods \cite{LevelSet, Cecil} and fast-marching methods \cite{Tsitsiklis,Sethian}, and have been used for optimal path planning problems in various applications including, among others, human walking paths \cite{Parkinson,Parkinson2,Arnold,Chen,Cartee2}, simple self-driving cars \cite{TakeiTsai1,TakeiTsai2,ParkinsonCar1,ParkinsonCar2}, planetary rovers \cite{Gee}, and pursuit-evasion games \cite{Cartee1}. These methods have the advantage that they are rooted in PDE and optimal control, from which one gains solid theoretical understanding of many facets of the algorithms, and thus robustness and interpretability are somewhat guaranteed. Classically, the disadvantage of these methods has been their relative lack of efficiency. Because they rely on grid-based methods for approximating solutions of PDE, they are not real time applicable even in low-dimensional problems, and suffer from the curse of dimensionality, making them entirely infeasible for high-dimensional problems. 

Recent numerical methods for Hamilton-Jacobi (HJ) type PDE attempt to break the curse of dimensionality using Hopf-Lax formulas \cite{Darbon,Lin}. Hamilton-Jacobi equations arise naturally from optimal control and dynamic programming, where they describe the behavior of a value function \cite{Bertsekas,Fleming,Liberzon}. Under mild conditions on the Hamiltonian $H:\R^d \to \R$ and the inital function $g:\R^d \to \R$, the classical Hopf-Lax formula \cite[Chapter 3.3]{EvansPDE} gives the solution of the HJ equation $
u_t + H(\nabla u) = 0,$ with initial condition $u(x,0) = g(x)$, in terms of a minimization problem. Exploiting this, one can resolve the solution of \eqref{eq:HJ} at individual points by solving an optimization problem, thus averting the need for discretized, grid-based approximations. However, this formula applies only to Hamiltonians which are state-independent, disallowing cases where $H$ depends explicitly on $x$, which account for essentially all interesting models of motion. The authors of \cite{Chow} provide a conjectured Hopf-Lax type formula for the state-dependent HJ equation \begin{equation} \label{eq:HJ} u_t + H(x,\nabla u) = 0, \,\,\,\,\,\, u(x,0) = g(x), \end{equation} where $H:\R^d \times \R^d \to \R$. We write $H = H(x,p)$ using $p$ as a proxy for $\nabla u$. Defining \begin{equation} \label{eq:Lagrangian}
L(x,p) = \langle p, \nabla_p H(x,p)\rangle - H(x,p),
\end{equation} they provide solid empirical and numerical evidence (and prove under restricted assumptions) that the solution of \eqref{eq:HJ} is given by \begin{equation} \label{eq:HLForm}
u(x,t) = \inf_{y \in \R^d}\left\{g(\B x(0)) + \int^t_0 L(\B x(s),\B p(s)) ds  \right\}
\end{equation} subject to the following Hamiltonian dynamics for $0 \le s < t$: \begin{equation} \label{eq:Hamdynamics}
\begin{split}
\dot{\B x}(s) &= \hphantom{-}\nabla_p H(\B x(s),\B p(s)), \\
\dot{\B p}(s) &= -\nabla_x H(\B x(s),\B p(s)), \\
\B x(t) &= x, \,\,\,\,\, \B p(t) = y.
\end{split}
\end{equation} 

Using this formula, we present a minimal-time optimal path planning method which is efficient enough to be semi-real-time applicable, while maintaining the PDE and optimal control formulation. To the authors' knowledge, this is the first PDE-based optimal path planning method which is applicable in nearly real time. As such, the method is both efficient and interpretable, and represents a significant step in averting the tradeoff between these two desirable properties. We develop the method with specific application to isotropic motion in two-dimensions, though with slight tweaks it could be applied more broadly to include more involved models of motion and more realistic modeling concerns.  

\section{MINIMAL-TIME OPTIMAL CONTROL}

In this section we give a brief inroduction to the minimal-time optimal control problem, and a formal derivation of the Hamilton-Jacobi-Bellman equation. Similar discussions with varying level of rigor are carried out in several texts \cite{Fleming,Bertsekas,Liberzon}.

Given a starting point $x \in \R^d$ and a desired final point $x_f \in \R^d$, the basic goal of minimal-time optimal control is to steer a controlled trajectory $\B x(\cdot)$ from $x$ to $x_f$ in the least possible time. Specifically, we assume that $\B x(\cdot)$ obeys \begin{equation}
\label{eq:dynamics}
\dot{\B x} = f(\B x,\Ba)\1_{x \neq x_f}(\B x)
\end{equation} for some dynamics function $f: \R^d \times A \to \R^d$. Here $A \subset \R^m$ is the set of admissible control actions, and $\B a(\cdot)$ is the control map taking values in $A$. It is assumed that the control map is chosen by some external user, and may be thought of as the ``steering plan" for the trajectory. Lastly, $\1_{x \neq x_f }$ is an indicator function taking value $1$ if $x\neq x_f$ and value $0$ if $x = x_f$. The inclusion of this indicator function is somewhat artificial, but it guarantees that once the trajectory reaches the final point, it stops moving, and as we will see shortly, it simplifies computations.

Given a horizon time $T$, the cost functional being minimized for the first arrival time problem is \begin{equation}
\label{eq:cost}
\mathcal C[\B x(\cdot),\B a(\cdot)] = \iota_{x_f}(\B x(t)) + \int^T_0 \1_{x \neq x_f}(\B x(s))ds,
\end{equation} where $\iota_{x_f}$ is the convex indicator function of the final point, taking value $+\infty$ if $x \neq x_f$ and value $0$ if $x = x_f$. Thus any path which does not reach the desired final point by time $T$ is assigned cost $+\infty$, and will never be optimal. Assuming a path does reach the final point by time $T$, the cost will be the first arrival time, since the indicator function in the integral will "turn off" once the path reaches $x_f$.  

We define the value function \begin{equation} 
\label{eq:value}
u(x,t) = \inf_{\B a(\cdot)} \mathcal C_{x,t}[\B x(\cdot),\B a(\cdot)]
\end{equation} where $\mathcal C_{x,t}$ is the same cost functional, restricted to the time interval $[t,T]$ and trajectories satisfying $\B x(t) = x$. Here $u(x,t)$ represents the remaining travel time required for a trajectory which is at $x$ at time $t$ to reach $x_f$ (assuming that there exists such a trajectory which can reach $x_f$ before time $T$; otherwise $u(x,t) = +\infty$). The standard dynamic programming argument, originally due to Bellman \cite{Bellman}, shows that $u(x,t)$ formally satisfies the terminal-valued Hamilton-Jacobi-Bellman equation \begin{equation}\label{eq:preHJB} \begin{split} &u_t + \inf_{a \in A}\{\1_{x \neq x_f}(\B x(s)) \langle f(x,a), \nabla u\rangle + \1_{x \neq x_f}(\B x(s)) \} = 0, \\ 
&u(x,T) = \iota_{x_f}(x).
\end{split}\end{equation} For a rigourous derivation, including a discussion of viscosity solutions, see \cite{Fleming}. 

Because initial-value problems are more familiar, we make the substitution $t \mapsto T-t$ (though in an abuse of notation, we still call the reversed-time value function $u$), and define the Hamiltonian \begin{equation}
\label{eq:Hamiltonian}
H(x,p) = -\1_{x\neq x_f}(x)\inf_{a \in A}\{\langle -f(x,a),p\rangle + 1\},
\end{equation} to arrive at \begin{equation}
\label{eq:HJB}
\begin{split}
&u_t + H(x,\nabla u) = 0, \\ 
&u(x,0) = \iota_{x_f}(x).
\end{split}
\end{equation} Assuming the solution of \eqref{eq:HJB} is known, the optimal feedback control $\B a^*(x,t)$ is defined as the argument achieving the minimum in \eqref{eq:Hamiltonian} (when $p$ is replaced by $\nabla u(x,t)$) and one can synthesize the optimal path by integrating \eqref{eq:dynamics} using the optimal feedback control. The solution to this equation has the somewhat special property that for each fixed $x$, $u(x,t)$ will become constant in finite time. Indeed, if the optimal travel time from $x$ to $x_f$ is $t^*$, then $u(x,t) = t^*$ for any $t \ge t^*$, and given time $t > t^*$ to travel, the optimal path will arrive at $x_f$ at time $t^*$, and then sit still. When solving, this is convenient because it means that the time horizon chosen for the problem is essentially arbitrary; it simply needs to be large enough that there is a path which connects $x$ and $x_f$ in the alloted time. This distinguishes the work in the manuscript from that in \cite{ParkinsonCar3}, where similar scalable numerical methods are developed but, because of a different modeling philosophy, one must choose the time horizon to be very near to the actual optimal travel time in order to resolve an approximation of an optimal path. While the methods presented in \cite{ParkinsonCar3} are as efficient and interpretable as those used here, needing to know the optimal travel time in advance precludes them from being applied in any real-time scenarios.  

We note that \eqref{eq:HJB} should only hold for $x \neq x_f$, and should be appended by a boundary condition $u(x_f,t) = 0$ for all $t$, as seen in \cite{ParkinsonCar2,Cartee1}. In that formulation, the $\1_{x \neq x_f}(\B x(s))$ in \eqref{eq:Hamiltonian} is superfluous. However, because we will numerically solve \eqref{eq:HJB} by recasting it as an optimization problem, it is convenient to avoid boundary conditions, which become constraints in the optimization problem. This is one reason to include the indicator function in \eqref{eq:dynamics}. Another is that in many examples, such as that of the simple self-driving car as in \cite{TakeiTsai1,TakeiTsai2,ParkinsonCar1} or isotropic motion as in \cite{Gee,Cartee1} one can explicitly resolve the optimal control values from the value function, and it is convenient to avert the need for special considerations at $x_f$. Using isotropic motion in $\R^d$ as an example, assume that $\B x(\cdot)$ may travel in any direction $a \in \mathbb S^{d-1}$, but with speed bounded by some function $v(x)>0$ defined throughout the domain. Then the dynamics are $\dot{\B x} = v(\B x)\B a$ where $\vert\B a(\cdot)\vert \le 1$. In this case, \begin{equation*}
-\inf_{\vert a \vert \le 1} \{v(x) \langle -a,\nabla u\rangle + 1\} = v(x)\vert \nabla u\vert - 1
\end{equation*} whenever $\nabla u$ exists, and the Hamiltonian is given by \begin{equation} \label{eq:IsoHamiltonian}
H(x,p) = \1_{x\neq x_f}(x)(v(x) \vert p \vert -1).
\end{equation} When $\nabla u$ exists, the optimal control is given by $a = -\nabla u/ \vert \nabla u\vert$. However, viscosity solutions of HJ equations may have points of non-differentiability, and in our case, this will occur at $x = x_f$ (perhaps among other locations). At $x = x_f$, one should take $a = 0$, so as to halt movement, but to avert the need for this special consideration, we can use the indicator function as in \eqref{eq:dynamics}. Having done so, the Hamilton-Jacobi-Bellman equation for isotropic motion (which is the equation we use in all examples below) is \begin{equation}
\label{eq:IsoHJB} \begin{split}
&u_t + \1_{x \neq x_f}(x)(v(x)\vert \nabla u \vert -1) = 0, \\
&u(x,0) = \iota_{x_f}(x). \end{split}
\end{equation} The remainder of this manuscript is concerned with approximating \eqref{eq:HJB} (with special application to \eqref{eq:IsoHJB}) using non-grid based numerical methods similar to those in \cite{Darbon,Lin}. We note that these methods can be applied more broadly, but as presented, it is important that the infimum in the Hamiltonian \eqref{eq:Hamiltonian} can be resolved so that the Hamiltonian is an explicit function of $x$ and $p$. As seen in \eqref{eq:IsoHamiltonian}, this is possible for the example of isotropic motion. It is also possible in any bang-bang control problem, including optimal trajectories for curvature constrained motion which can be used to model simple vehicles \cite{TakeiTsai2,ParkinsonCar3}. An example of an application where this is not possible is the model for human walking paths in mountainous terrain in \cite{Parkinson,Parkinson2}, where the speed of motion is assumed to depend on local slope of the terrain. Applying this method to situations like that would require more work. 

\section{NUMERICAL METHODS}

In this section, we provide a brief exposition of the numerical methods used to solve \eqref{eq:HJ} using the Hopf-Lax formula described in \eqref{eq:HLForm} and \eqref{eq:dynamics}. We present an algorithm that computes optimal paths very efficiently, and can include semi-real-time adjustment and recalculation, meaning that it can account for real time concerns like obstacle discovery. Such problems are notoriously difficult to tackle using the feedback control / Hamilton-Jacobi formulation, because any new information (for example, discovery of an obstacle) requires a new PDE solve. To the authors' knowledge, this is the first algorithm which maintains the Hamilton-Jacobi formulation, but is applicable to real-time problems.

We solve the minimization problem described in \eqref{eq:HLForm} and \eqref{eq:dynamics} using the splitting method described in \cite{Lin}, which in turn employs a primal-dual algorithm in the spirit of \cite{PDHG}.  The construction of the splitting method is coved in detail in \cite{Lin}. We provide a brief overview and discuss its specific application to our problem. 

The splitting method works by discretizing path-space, and alternately minimizing the Hamiltonian with respect to the state variables and the co-state variables (which are proxies for $\nabla u$ along the path), and iterating until convergence. Specifically, to approximate the solution of \eqref{eq:HJ} and a point $(x,t) \in \R^{d}\times[0,\infty)$, we first discretize the interval $[0,t]$ into $J$ smaller intervals of length $\delta = t/J$. Let $t_j= j\delta$, for $j = 0,1,\ldots,N$, and let $x_j$ and $p_j$ for $j = 0,1,\ldots,N$ be approximations to the points $\B x(t_j)$ and $\B p(t_j)$ along the path. Beginning from a Langrangian formulation and working formally, \cite{Lin} derives an approximation of the solution of \eqref{eq:HLForm} and \eqref{eq:dynamics} (and thus the solution of \eqref{eq:HJ}) in the form of a saddle point problem: \begin{equation}\label{eq:saddleProb}
u(x,t) \approx \max_{\{p_j\}} \min_{\{x_j\}} \left\{ 
\begin{split} g(x_0) + &\sum^J_{j=1}\innerprod{p_j}{x_j - x_{j-1}} \\ &\hspace{1cm}- \delta\sum^J_{j=1} H(x_j,p_j)\end{split}\right\}.
\end{equation} 

This optimization problem is solved using alternating primal-dual optimization as described in algorithm \ref{alg:1}, which is adapted from \cite{Lin}. The algorithm takes advantage of the fact that the minimization over each individual vector $x_j$ or $p_j$ in \eqref{eq:saddleProb} is decoupled from the others, so the formula can be minimized with respect to each vector along the path individually, rather than with respect to the entire path at once. One of the great strengths of this method is that, while resolving the saddle point problem \eqref{eq:saddleProb} gives the value of of the solution of \eqref{eq:HJ} (in our case the optimal travel time), $\{x_j\}_{j=0}^{J}$ represents an approximation of the optimal path. 

We now describe the specific application of this method to \eqref{eq:IsoHJB}. In this case, we are using the Hamiltonian $H(x,p)$ given by \eqref{eq:IsoHamiltonian}, and the inital function $g(x) = \iota_{x_f}(x)$, the convex indicator of the desired final point. Note that $g$ only appears in algorithm \ref{alg:1} in the update for $x_0^{k+1}$. Because the convex indicator takes value $+\infty$ at any point which is not $x_f$, the minimum can only occur at $x_f$, so we will always have $x_0^{k+1} = x_f$. Since there are roughly $2J$ optimization problems at each iteration, solving each optimization problem efficiently is of utmost importance. As it turns out, the update for $p^{k+1}_j$ can actually be resolved explcitly: \begin{equation} \label{eq:pupdate}
p^{k+1}_j = \max\left\{0, 1 - \frac{\sigma \delta \1_{x\neq x_f}(x_j^k) v(x_j^k)}{\vert \beta^k_j\vert}\right\} \beta^k_j,
\end{equation} where $\beta_j^k$ is as in algorithm \ref{alg:1}. This formula is derived in \cite[\S3.1]{ParkinsonCar3}, in a slightly different context, though the details are essentially the same. The optimization problem for $x^{k+1}_j$ cannot be resolved explicitly except in the case of very simple velocity functions $v(x)$. In this case, one needs to approximately solve \begin{equation} \label{eq:xupdate} x^{k+1}_j =  \argmin_{\tilde x}\{ -\delta\tau H(\tilde x,p_j^{k+1}) + \frac{1}{2} \lvert \tilde x -\nu_j^k \rvert_2^2 \}.\end{equation} Empirically, \cite{Lin} found that this approximation could be quite crude, and the algorithm still works. This was corroborated by our implementation where, to update $x^{k+1}_{j}$, we start from $x^k_{j}$ and simply take one gradient descent step with rate $\gamma$. Thus the update we use is \begin{equation} \label{eq:xupdate2} x^{k+1}_j = x^k_j - \gamma(-\delta \tau \nabla_x H(x_j^k,p^{k+1}_j )+ (x^k_j-\nu^k_j)). \end{equation} Note that $p^{k+1}_j$ has already been resolved when we arrive at this step, so this is entirely explicit. 

\begin{algorithm}[t!]
\caption{Splitting Method for Solving \eqref{eq:saddleProb}}
Input the point $(x,t)$ at which to resolve the HJ equation, as well as the max iteration count $k_{\text{max}}$, proximal step sizes $\sigma, \tau$, relaxation parameter $\kappa$, and convergence tolerance TOL. \\

 Set $x^1_J= x, p^1_0 = 0$. Initialize $\{x^1_j\}^{J-1}_{j=0}$, $\{p^1_j\}^J_{j=1}$ randomly, and set $z^1_j = x^1_j$ for all $j = 0,1,\ldots, J$. 
\begin{algorithmic}[t!]
\For {$k = 1$ to $k_{\text{max}}$}\\
    \State $p^{k+1}_0 = 0$
    \For {$j = 1$ to $N$}
    
    \State $\beta_{j}^k = p_j^k + \sigma(z_j^k - z_{j-1}^k)$
    \State $p_j^{k+1} = \argmin_{\tilde p} \{ \delta H(x_j^k,\tilde p) + \frac{1}{2\sigma} \lvert \tilde p- \beta^k_j \rvert^2 \}$
    \EndFor\\
    
    \State $\nu_{0}^k = x_0^k + \tau p_1^{k+1}$
    \State $x_0^{k+1} = \argmin_{\tilde x} \{ g(\tilde x) + \frac{1}{2\sigma} \lvert \tilde x - \nu^k_{0} \rvert^2 \}$
    \For{$j = 1$ to $N - 1$}
    
    \State $\nu_j^k = x_j^k - \tau(p_j^{k+1} - p_{j+1}^{k+1})$
    \State $x_j^{k+1} = \argmin_{\tilde x}\{ - \delta H(\tilde x,p_j^{k+1}) + \frac{1}{2\tau} \lvert \tilde x -\nu_j^k \rvert^2 \}$
    \EndFor
    \State $x^{k+1}_N = x$ \\
    
    \For{$j = 0$ to $N$}
    \State $z_j^{k+1} = x_j^{k+1} + \kappa(x_j^{k+1} - x_j^k)$
    \EndFor\\
    
    \State change $= \max \{ \lVert x^{k+1} - x^k \rVert , \lVert p^{k+1} - p^k \rVert \} $
    \If {change $<$ TOL}
    \State break
    \EndIf
\EndFor
\State $u = g(x_0) + \sum_{j=1}^N \langle p_j, x_j - x_{j-1} \rangle - \delta H(x_j,p_j)$
\State \textbf{return } $u$, $\{x_j\}_{j=0}^{J}$
\end{algorithmic}
\label{alg:1}
\end{algorithm}

\subsection{ACCOUNTING FOR OBSTACLES AND \\APPROXIMATING INDICATOR FUNCTIONS} 

One final concern is how to include impassible obstacles in the model. A common method of doing this when solving motion planning problems in the Hamilton-Jacobi formulation is to simply set the value function to $+\infty$ inside obstacles. This is akin to assigning infinite cost to any path which intersects and obstacle, so that such a path would never be optimal. However, because we are not solving in a grid-based manner, whereupon the value function $u(x,t)$ is resolved from nearby points, this manner of including obstacles is not available to us. One possibility would be to simply restrict the domain for the argmin in \eqref{eq:xupdate}. However, this would lead to a difficult state-constrained optimization problem. Instead, we deal with obstacles in a manner similar to \cite{TakeiTsai2,ParkinsonCar3}: we set the velocity to zero inside obstacles. To do so, we multiply $f(x,a)$ in equation \eqref{eq:Hamiltonian} or $v(x)$ in equation \eqref{eq:IsoHamiltonian} by a function $O(x)$ which takes value 1 in the free space and value 0 inside obstacles. This follows through all the calculations in an entirely predictable manner; for example, an $O(x^k_j)$ will simply appear in the numerator in \eqref{eq:pupdate}.

The one place where this requires some additional thought is in the update rule for $x^{k+1}_j$ given by \eqref{eq:xupdate2}, where $\nabla_xH(x,p)$ appears. Having inserted the discontinuous function $O(x)$ into \eqref{eq:IsoHamiltonian}, $\nabla_xH(x,p)$ is no longer well-defined at the boundarys of obstacles. In fact, because of the presence of $\1_{x\neq x_f}(x)$, $H(x,p)$ is already non-differentiable at $x = x_f$. This latter problem can be handled by explicitly checking if $x_f$ is the minimizer of \eqref{eq:xupdate}, and only performing the gradient descent \eqref{eq:xupdate2} if the minimizer lies elsewhere. Empirically, having tried this, the number of iterations required for convergence was exorbitantly large, leading to inefficiency. Instead, we approximate both the indicator function of the obstables $O(x)$ and the indicator function $\1_{x\neq x_f}(x)$ by smooth functions. For the indicator of the final point, we use \begin{equation} \label{eq:indApprox} 
\1_{x\neq x_f} (x) \approx 1 - e^{-A\vert x - x_f \vert^2}
\end{equation} for some large parameter $A$.  For the obstacles, we use \begin{equation} \label{eq:obsApprox}
O(x) \approx \frac 1 2 + \frac 1 2 \tanh(B d(x))
\end{equation} where $d(x)$ is the signed distance to the boundary of the obstacles (negative inside the obstacles), and $B$ is another large parameter. In all of our examples, we take $A = B = 100$. This introduces one final difficulty of how to efficiently compute $d(x)$, which can itself be expressed as the solution of a Hamilton-Jacobi equation. To deal with this, we use the strategy of \cite{ParkinsonCar3}, where all obstacles are required to be circles, or to be approximated by a collection of disjoint circles. This allows for very efficient calculation of both $d(x)$ and $\nabla d(x)$ (which will appear in $\nabla_x H(x,p)$) when $x$ is outside obstacles as described in \cite{ParkinsonCar3}. For the moment, we model stationary obstacles only; investigation into the application of these methods to scenarios with moving obstacles is ongoing. 

Because of these smooth approximations to the indicator functions, the gradient $\nabla_xH(x,p)$ is large near $x_f$ and near obstacle boundaries. Accoringly, the gradient descent rate $\gamma$ in \eqref{eq:xupdate2} needs to be very small to resolve the path near these points. Away from these points, where it is easier to resolve the path, it is better to take a large gradient step in order to accelerate convergence. To account for this, one can either let the gradient descent rate depend on the spatial variable $x$, or begin with a large gradient descent rate to resolve the basic skeleton of the path, and then decrease the rate for higher iterations so as to refine the path at points which are difficult to resolve.

\subsection{REAL-TIME CORRECTION} 

It was mentioned briefly, though it bears repeating, that one of the biggest strengths of the modeling decisions we have made in this manuscript is that the time $t$ that one chooses is arbitrary, so long as it is large enough that the optimal path requires time less than $t$ to traverse. For example, if the optimal path from $x$ to $x_f$ requires 5 seconds to traverse, then inputting any of $t = 5, 8$ or $20$ will all provide the same result: a path which reaches $x_f$ by time $t = 5$, and then sits still until the time limit is reached. This distinguishes the work in this manuscript from that in \cite{ParkinsonCar3}, where similar methods are used, but due to different modeling decisions, one must know the optimal travel time in advance. Because it is not necessary to know the optimal travel time when computing the path, this method is amenable to real-time correction when new information presents itself. This could account for any number of realistic scenarios, such as obstacle discovery as demonstrated in examples below. In the examples, there are obstacles present, but the traveler does not know about them until they are within some fixed radius. Thus, at the beginning, we compute the optimal path as if no obstacles are present, then begin traveling, and recompute a path any time an obstacle is encountered. With real-time planning like this, we forsake the hope of globally optimal trajectories, since there is imperfect information, but at every juncture, the traveler computes an optimal path given the current information. Such methods could be implemented online in real world optimal routing applications. To the authors' knowledge, this is the first algorithm which maintains the Hamilton-Jacobi formulation, but is efficient enough (even for high-dimensional systems) to be implemented online.  

\section{RESULTS AND CONCLUSIONS}

We present some examples of the application of our algorithm to isotropic motion.  In all cases we operate in two spatial dimensions, set $t = 8$ and $\delta = 0.1$, the starting point $x = (-1,-1)$, and the final point $x_f = (1,1)$. The parameters for algorithm \ref{alg:1} are set at $\sigma = 1, \tau = 0.2, \kappa = 1$. The convergence tolerance is TOL $=10^{-3}$, and the maximum iteration count is 40000, though this is seldom reached (empirically, very few paths require more than 25000 iterations to resolve).  For our gradient descent rate in \eqref{eq:xupdate2}, we begin with $\gamma = 0.2$ for the first 5000 iterations, and halve $\gamma$ every 1000 iterations thereafter. For obstacle discovery, we halt and recompute the trajectory any time the vehicle is within distance $0.1$ of a new obstacle. In these examples the times and spatial distances are artificial. The real contribution of this work is not high-fidelity modeling of any one real world system, but rather proof of concept that one can indeed perform nearly real-time path planning in the Hamilton-Jacobi formulation. 

\begin{figure*}[!]
\fbox{\includegraphics[width=0.18\textwidth,trim = 60 40 50 30, clip]{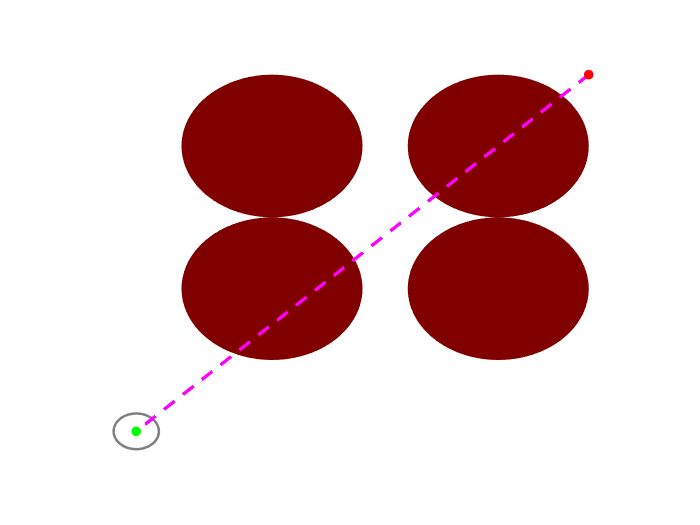}} \, 
\fbox{\includegraphics[width=0.18\textwidth,trim = 60 40 50 30, clip]{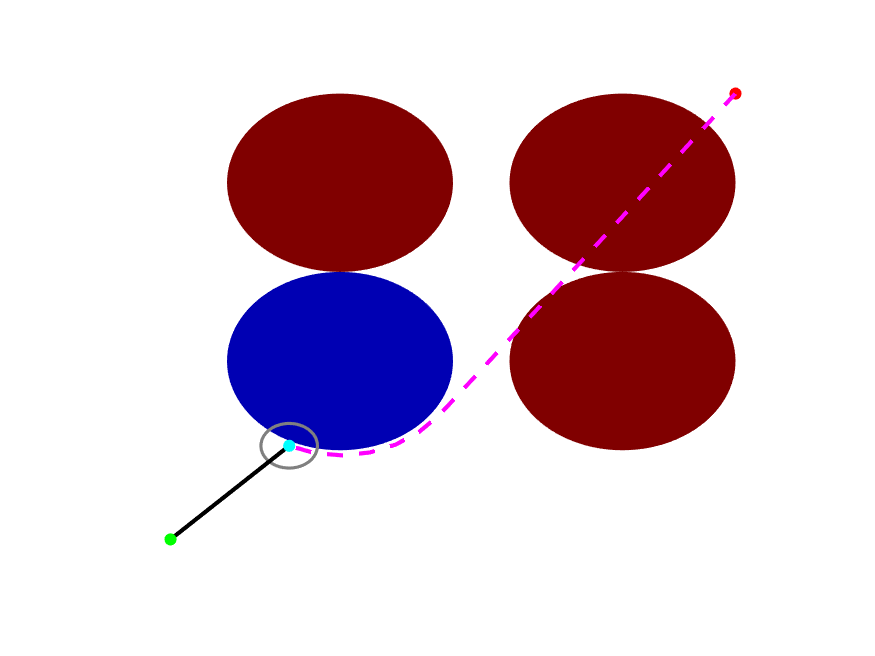}} \, 
\fbox{\includegraphics[width=0.18\textwidth,trim = 60 40 50 30, clip]{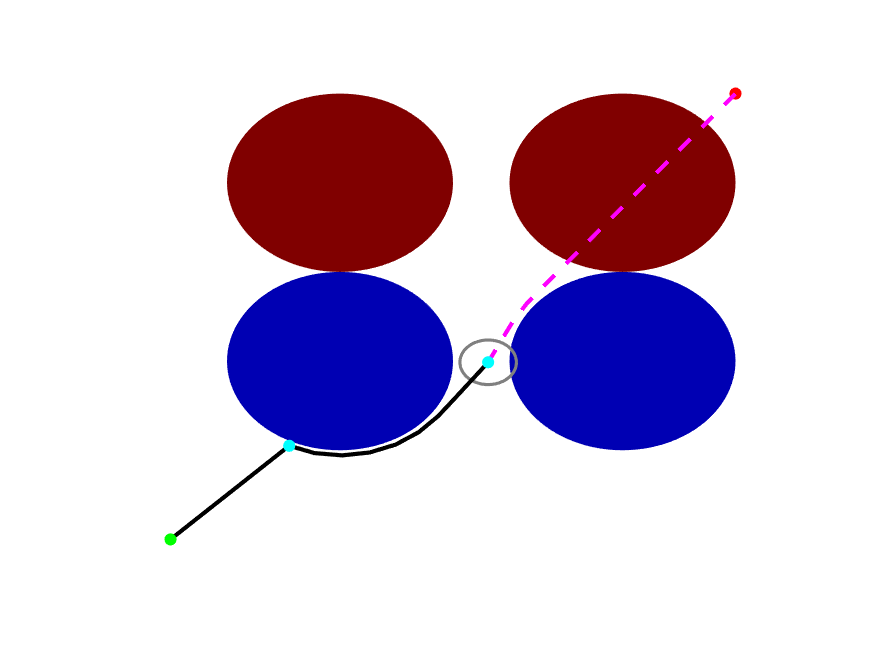}}\, 
\fbox{\includegraphics[width=0.18\textwidth,trim = 60 40 50 30, clip]{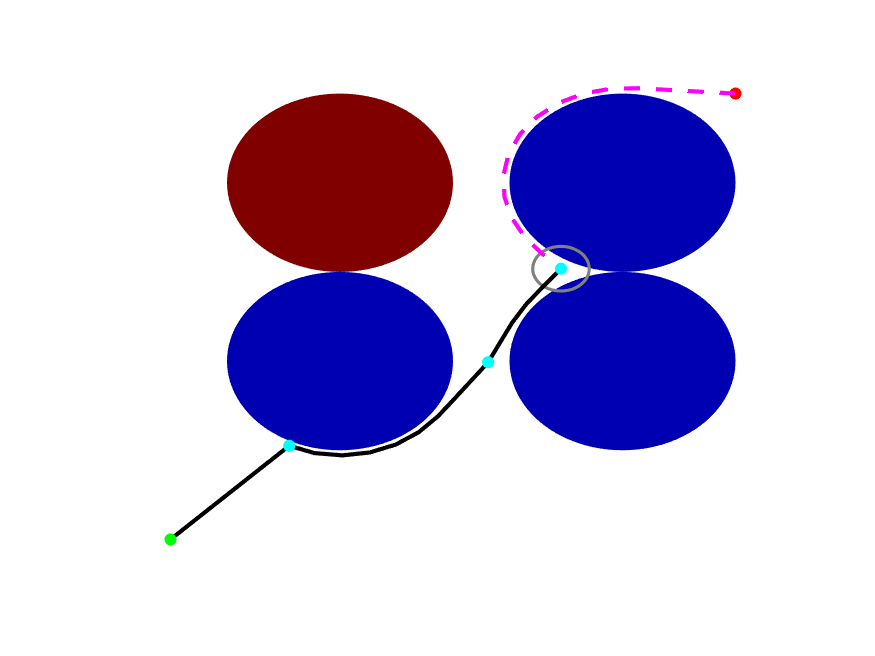}} \, 
\fbox{\includegraphics[width=0.18\textwidth,trim = 60 40 50 30, clip]{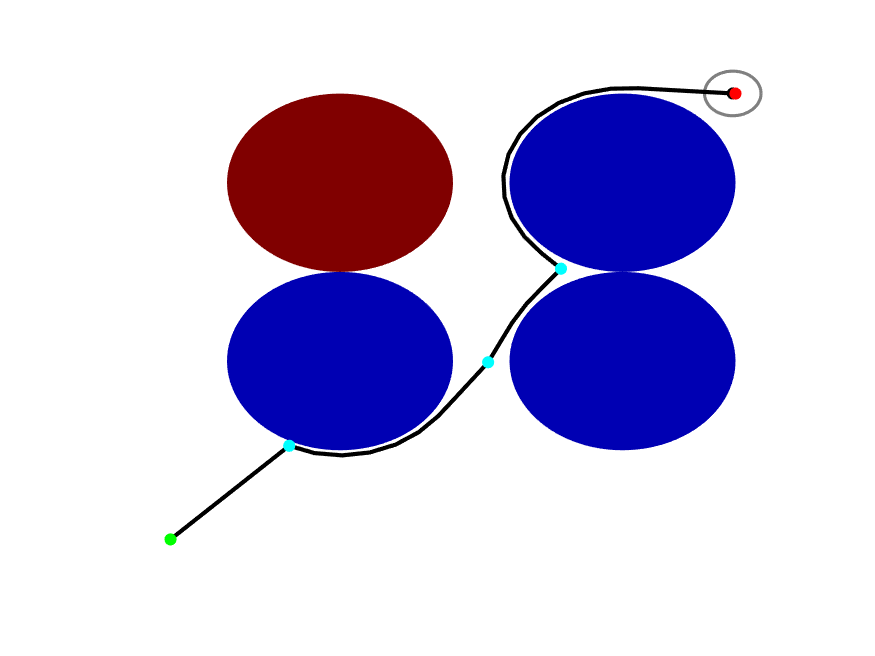}}
\caption{A vehicle traveling with constant velocity $v(x) = 1$ navigates around obstacles.}
\label{fig:1}
\end{figure*}

\begin{figure*}[!]
\fbox{\includegraphics[width=0.18\textwidth,trim = 60 40 50 30, clip]{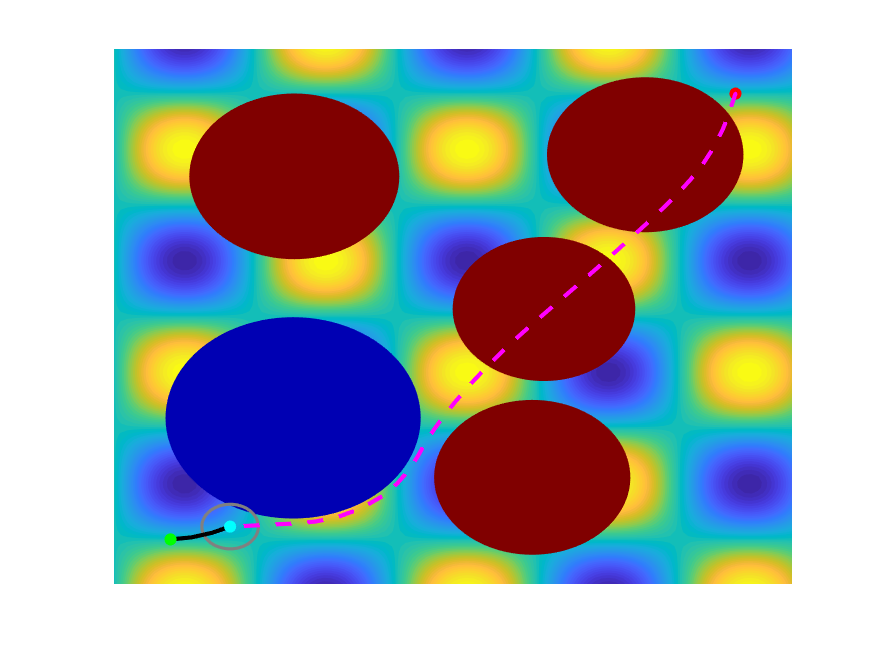}} \, 
\fbox{\includegraphics[width=0.18\textwidth,trim = 60 40 50 30, clip]{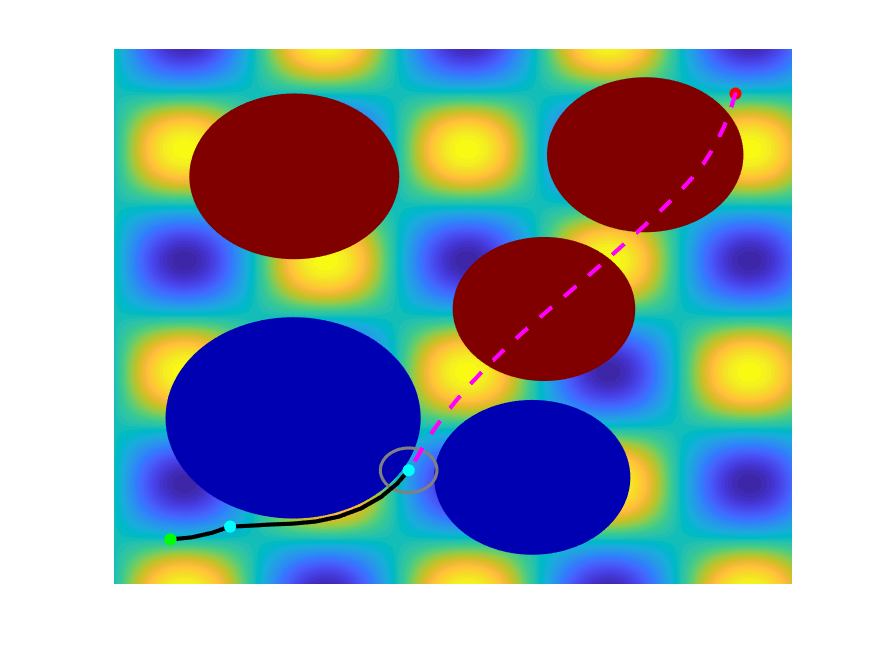}} \, 
\fbox{\includegraphics[width=0.18\textwidth,trim = 60 40 50 30, clip]{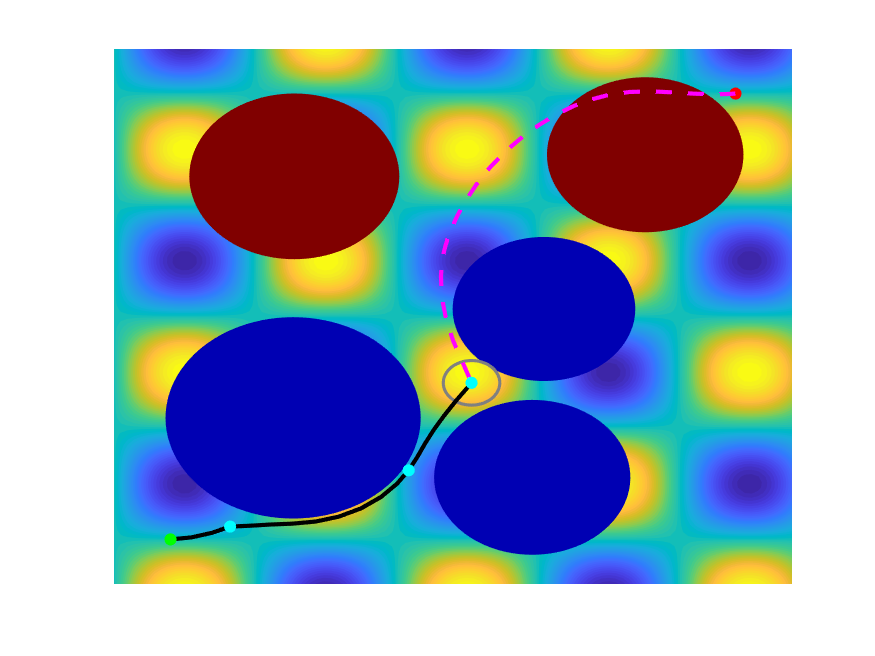}}\, 
\fbox{\includegraphics[width=0.18\textwidth,trim = 60 40 50 30, clip]{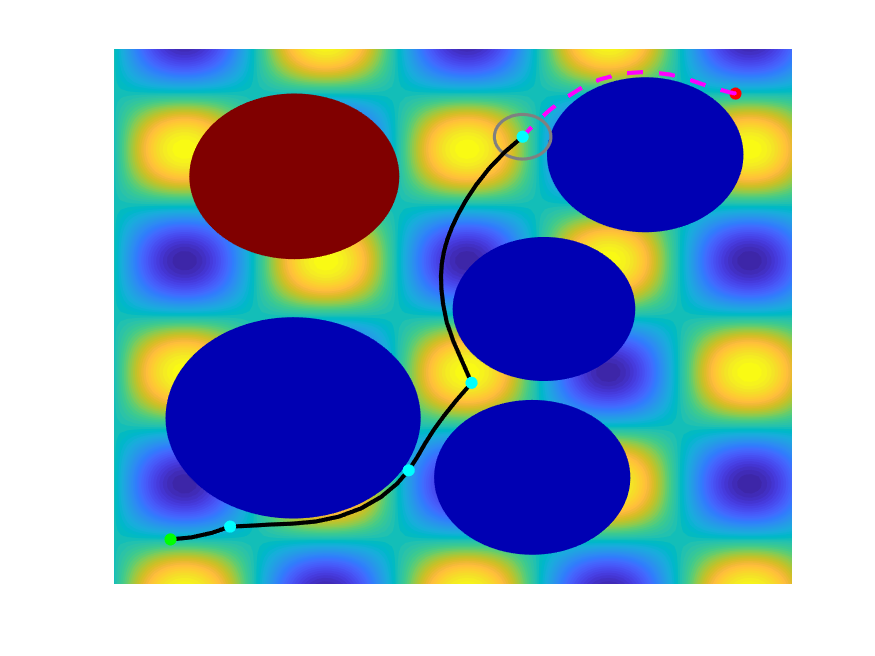}}\, 
\fbox{\includegraphics[width=0.18\textwidth,trim = 60 40 50 30, clip]{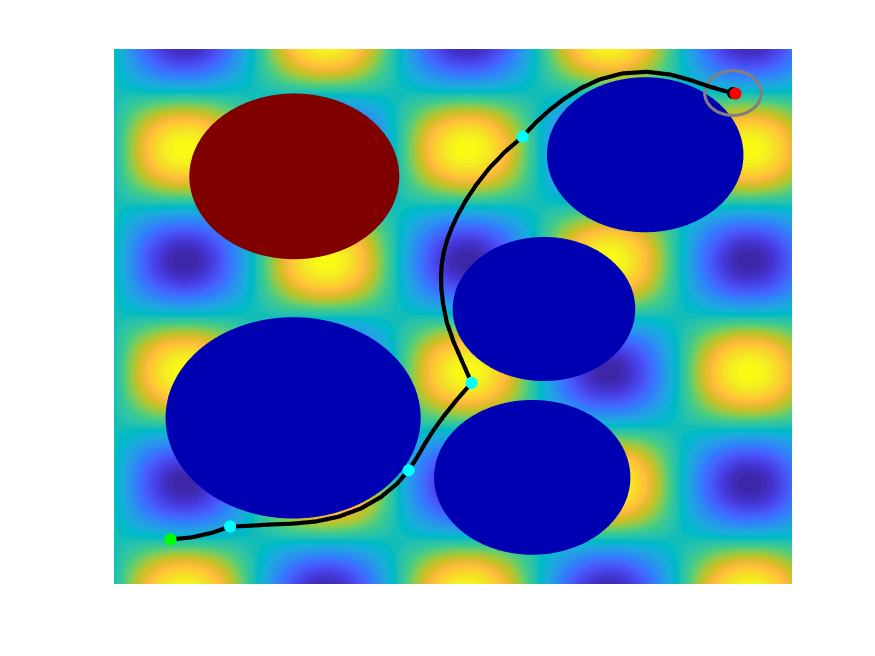}}
\caption{A vehicle traveling with velocity $v(x) =\tfrac 3 4 - \tfrac 1 4 \sin(2\pi (x_1+\tfrac 1 5)) \sin(2\pi x_2)$ navigates around obstacles. The velocity function is plotted in the background (yellow represents high velocity regions; blue represents low velocity regions).}
\label{fig:2}
\end{figure*}

We present two examples. In each plot, the vehicle is the cyan dot, and we leave a cyan dot behind at every point where an obstacle was discovered and a new path was computed. The starting point is the green dot and the ending point is the red dot. Obstacles are plotted in dark red if they are undiscovered, and blue once they are discovered. The black curve is the portion of the path which has already been traveled. The dotted magenta curve is the planned path given the current information. In the first example (seen in fig. \ref{fig:1}), the velocity function is the constant $v(x) = 1$; in the second example (seen in fig. \ref{fig:2}), it is the oscillatory function $v(x) = \tfrac 3 4 - \tfrac 1 4 \sin(2\pi (x_1+\tfrac 1 5)) \sin(2\pi x_2)$. The code which produced these is posted on GitHub, and displays these as animations which are preferable for demonstrating the algorithm.\footnote{\url{https://github.com/chparkin/SRT_HJ_Alg}}

In both cases, the vehicle successfully reaches the goal, after halting and recomputing the path several times when encountering obstacles. For the first example, each path was resolved by algorithm \ref{alg:1} in 10000 iterations or fewer. With the more complicated velocity function in the second example, roughly 25000 iterations were required to resolve each path. Due to the random initialization, there is some mild stochasiticity in these results, but in our simulations, no trajectory for either of the above examples ever failed to reach the final point, and no individual path failed to resolve within the maximum iteration count of 40000. On the first author's personal laptop (Intel(R) Core(TM) i7-10510U processor running at 1.80GHz, 12GM RAM), roughly 1 second of computation time was required for every 5000 iterations, so for the first example, each path was resolved in roughly 2 seconds, and for the second example, each path was resolved in roughly 5 seconds. We note that note that the convergence could likely be accelerated by resolving the minimization \eqref{eq:xupdate} more accurately, or performing a more rigorous study of the manner in which each parameter affects convergence, as suggested by \cite{Lin}, so it is very likely these computation times could be improved. Because updated paths could be computed online as a vehicle is moving, this algorithm is efficient enough for many real-time applications.  

In this manuscript, we present a general control-theoretic framework for minimum-time path planning, based on dynamic programming and a Hamilton-Jacobi formulation. We also design an algorithm for resolving optimal paths based on Hopf-Lax type formulas which is efficient enough to be real-time applicable. This represents a significant step toward routing algorithms which are efficient, scalable, and fully interpretable. In ongoing work, the authors are exploring the application of these methods to higher dimensional systems, and to other models of motion (for example, simple vehicles and path planning on manifolds). Another interesting direction would be to apply similar methods to scenarios where minimizing travel time is not the only goal. For example, one may consider terrain discovery, energy efficiency, passenger comfortability, or any number of other real world concerns.

\addtolength{\textheight}{0cm}

\bibliographystyle{ieeetr}
\bibliography{Biblio}

\end{document}